\documentclass[10pt]{amsart}
\usepackage{amsfonts}
\usepackage{amsmath}
\usepackage{amsthm}
\usepackage{amssymb}
\usepackage{latexsym}
\usepackage{multicol}
\usepackage{verbatim}
\usepackage{tabularx}

\advance\textwidth by 1.2in \advance\oddsidemargin by -.6in
\advance\evensidemargin by -.6in
\newtheorem*{cor}{Corollary}
\newtheorem*{lem}{Lemma}
\newtheorem*{prop}{Proposition}

\theoremstyle{definition}
\newtheorem*{defn}{Definition}
\theoremstyle{definition}
\newtheorem*{thm}{Theorem}

\newenvironment{pf}{\proof}{\endproof}
\newcounter{cnt}
\newenvironment{enumerit}{\begin{list}{{\hfill\rm(\roman{cnt})\hfill}}{%
\settowidth{\labelwidth}{{\rm(iv)}}\leftmargin=\labelwidth%
\advance\leftmargin by
\labelsep\rightmargin=0pt\usecounter{cnt}}}{\end{list}}

\theoremstyle{remark}

\numberwithin{equation}{section} 

\newtheorem{example}{Example}

\def\mb#1{\text{$\mathbf{#1}$}}

\def\wt{{\rm wt}}

\def\ord#1^#2{#1$^{\text{#2}}$}

\def\wt{{\rm wt}}
\def\sop_#1^#2{\text{\scriptsize $\bigoplus\limits_{#1}^{#2}$}}

\begin{document}

\newcommand{\thmref}[1]{Theorem~\ref{#1}}
\newcommand{\secref}[1]{Section~\ref{#1}}
\newcommand{\lemref}[1]{Lemma~\ref{#1}}
\newcommand{\propref}[1]{Proposition~\ref{#1}}
\newcommand{\corref}[1]{Corollary~\ref{#1}}
\newcommand{\remref}[1]{Remark~\ref{#1}}
\newcommand{\defref}[1]{Definition~\ref{#1}}
\newcommand{\er}[1]{(\ref{#1})}
\newcommand{\id}{\operatorname{id}}
\newcommand{\tensor}{\otimes}
\newcommand{\nc}{\newcommand}
\newcommand{\rnc}{\renewcommand}
\newcommand{\qbinom}[2]{\genfrac[]{0pt}0{#1}{#2}}
\nc{\cal}{\mathcal} \nc{\goth}{\mathfrak} \rnc{\bold}{\mathbf}
\renewcommand{\frak}{\mathfrak}
\newcommand{\desc}{\operatorname{desc}}
\newcommand{\Maj}{\operatorname{Maj}}
\renewcommand{\Bbb}{\mathbb}
\nc\bpi{{\mbox{\boldmath $\pi$}}} \nc\bvpi{{\mbox{\boldmath
$\varpi$}}}
 \nc\balpha{{\mbox{\boldmath $\alpha$}}}
\newcommand{\lie}[1]{\mathfrak{#1}}
\makeatletter
\def\section{\def\@secnumfont{\mdseries}\@startsection{section}{1}%
  \z@{.7\linespacing\@plus\linespacing}{.5\linespacing}%
  {\normalfont\scshape\centering}}
\def\subsection{\def\@secnumfont{\bfseries}\@startsection{subsection}{2}%
  {\parindent}{.5\linespacing\@plus.7\linespacing}{-.5em}%
  {\normalfont\bfseries}}
\makeatother
\def\subl#1{\subsection{}\label{#1}}

\nc{\Cal}{\cal} \nc{\Xp}[1]{X^+(#1)} \nc{\Xm}[1]{X^-(#1)}
\nc{\on}{\operatorname} \nc{\ch}{\mbox{ch}} \nc{\Z}{{\bold Z}}
\nc{\J}{{\cal J}} \nc{\C}{{\bold C}} \nc{\Q}{{\bold Q}}
\renewcommand{\P}{{\cal P}}
\nc{\N}{{\Bbb N}} \nc\boa{\bold a} \nc\bob{\bold b} \nc\boc{\bold
c} \nc\bod{\bold d} \nc\boe{\bold e} \nc\bof{\bold f}
\nc\bog{\bold g} \nc\boh{\bold h} \nc\boi{\bold i} \nc\boj{\bold
j} \nc\bok{\bold k} \nc\bol{\bold l} \nc\bom{\bold m}
\nc\bon{\bold n} \nc\boo{\bold o} \nc\bop{\bold p} \nc\boq{\bold
q} \nc\bor{\bold r} \nc\bos{\bold s} \nc\bou{\bold u}
\nc\bov{\bold v} \nc\bow{\bold w} \nc\boz{\bold z}

\nc\ba{\bold A} \nc\bb{\bold B} \nc\bc{\bold C} \nc\bd{\bold D}
\nc\be{\bold E}  \nc\bg{\bold G} \nc\bh{\bold H}
\nc\bi{\bold I} \nc\bj{\bold J} \nc\bk{\bold K} \nc\bl{\bold L}
\nc\bm{\bold M} \nc\bn{\bold N} \nc\bo{\bold O} \nc\bp{\bold P}
\nc\bq{\bold Q} \nc\br{\bold R} \nc\bs{\bold S} \nc\bt{\bold T}
\nc\bu{\bold U} \nc\bv{\bold V} \nc\bw{\bold W} \nc\bz{\bold Z}
\nc\bx{\bold X}

\title{Branched Crystals and Category $\mathcal{O}$}
\author{Vyjayanthi Chari, Dijana Jakeli\'c, and Adriano A. Moura }
\address{Department of Mathematics, University of
California, Riverside, CA 92521} \email{chari@math.ucr.edu,
jakelic@math.ucr.edu, adrianoam@math.ucr.edu}
\address{Max-Planck-Institut f\"ur Mathematik, D-53111, Bonn, Germany} \email{jakelic@mpim-bonn.mpg.de} \maketitle
\setcounter{section}{0}

\flushbottom
\section*{Introduction}
The theory of crystal bases  introduced by Kashiwara in \cite{K}
to study the category of integrable representations of quantized
Kac--Moody Lie algebras has been a major
 development in the combinatorial
approach to representation theory. In particular Kashiwara defined
the tensor product of crystal bases and showed that it
corresponded to the tensor product of representations. Later, in
\cite{Ka} he  defined the abstract notion of a crystal, the tensor
product of crystals and showed that the tensor product was
commutative and associative.

 In this paper, using some ideas from \cite{J1} and \cite{Ze},
 we give a definition of branched crystals
adapted to the study  of the Bernstein--Gelfand--Gelfand category
$\cal O$ of the quantized enveloping algebra of $\lie{sl_2}$ which
coincides with the usual definition of crystals for the integrable
modules. In \cite{K},  Kashiwara  essentially  defined a crystal
basis for the Verma modules. However, it is not hard to see that
his tensor product rule, even for the case of $\lie{sl}_2$, does
not give the decomposition (as  a direct sum of indecomposable
modules) of the corresponding tensor product of modules. Also, the
restriction of the basis to a particular color does not reflect
the decomposition of the module as a direct sum of indecomposable
modules for the corresponding $\lie{sl}_2$. We define a notion of
the tensor product of  branched crystals which extends Kashiwara's
definition and the definition in \cite{Ze} and prove in a
combinatorial manner that the tensor product decomposes in the
same way as the corresponding representations. Using this, we are
then able to prove that the tensor product is both associative and
commutative. It is not hard to generalize the definition of
branched crystals to the higher rank case although the connection
with $\cal O$ is probably much harder to establish. However, the
results of Section 1 make it plausible that this is essentially
the only possible theory of crystals for $\cal O$ which would
satisfy the requirement that the crystal corresponding to a
representation $V$, when restricted to a particular color, is the
crystal of $V$ regarded as a module for $\bu_q(\lie{sl}_2)$.

 The paper is organized as follows. In Section 1 we define the
 notion of a (one--colored)  branched crystal, its indecomposable components  and
 classify the indecomposable branched crystals. In Section 2 we
 define the  tensor product of two branched crystals and show that
 the result is also a branched crystal.
Finally, we identify the indecomposable components (with
multiplicities) in the tensor product of indecomposable branched
crystals. In Section 3 we make a connection with the
representation theory of $\bu_q(sl_2)$. Using the results of
Section 1 we see that given a module $V$ in $\cal O$, we can
associate to it in a purely formal but natural and unique way a
branched crystal $B(V)$ so that direct sums are preserved.
Further, the results of Section 2 make it clear that this
association preserves tensor products, i.e., $B(V\otimes W)\cong
B(V)\otimes B(W)$ for all $V,W\in \cal {O}$.

{\bf Acknowledgements:} We thank J. Greenstein for bringing
\cite{Ze} to our attention and for useful discussions.

%\newpage
\section{Branched Crystals}\label{s:c}

In this section we introduce the notion of a (one colored)
branched crystal. These combinatorial objects are
 analogous to the notion  of normal  crystals defined by Kashiwara
in \cite{Ka}.

\subl{}

\label{bc}
\begin{defn}\label{crystal}
A branched crystal is a nonempty set  $B$
 together with maps
 $\tilde e,\tilde f:B\sqcup\{0\}\to B\sqcup\{0\}$, $\wt:B\to \bz$, $\varepsilon,\varphi:B\to \bz$
   satisfying the following axioms:
\begin{enumerate}
\item[(i)] $\tilde e0=\tilde f0=0$ and if $\tilde e b \in B$
(resp.$\tilde f b\in B$),  then $\wt(\tilde e b) = \wt(b)+2$
(resp.  $\wt(\tilde f b) = \wt(b)-2$)),
\item[(ii)] if $\tilde f b\in B$, then $\tilde e\tilde f b=b$,
\item[(iii)] $\tilde eb\ne \tilde eb'$ if $b,b'\in B$ are such that
$\tilde e b,\tilde eb'\in B$ but $\tilde f\tilde eb\ne b$ and
$\tilde f\tilde eb'\ne b'$,
\item[(iv)] $\varepsilon(b)=\min\{s\in \bz^+: \tilde f\tilde e ^{s+1}b\ne \tilde
e^sb\}$,
\item[(v)] $\varphi(b)=\wt(b)+\varepsilon (b)$. If there exists $l>0$ such that $\tilde
f^{l-1}b\ne 0$ and $\tilde f^lb=0$, then $l=\varphi(b)+1$. Further, if $b\in B$ is
such that  $\tilde eb\in B$ but $\tilde f\tilde eb\ne b$, then
$\varphi(\tilde eb)=0$.
%\item[(vi)] there exists $r\equiv r(B)\in\bz$ such that $\wt(b)\le r$ for all $b\in B$.
\end{enumerate}
\end{defn}

 An element $b\in B$ is called a branch point if $\tilde e b\in B$
but $\tilde f\tilde e b\ne b$ and  $B^{br}$ denotes the set of
branch points in $B$. Note that   $\tilde fb\notin B^{br}$
 for $b\in B$.
 A subcrystal of $B$ is a subset $B'$ of $B$
such that $B'$ together with the restrictions of the maps
$\wt,\tilde e,\tilde f,\varepsilon,\varphi$ is  a branched crystal. Given $b\in B$, let $B_b$ be the subcrystal $$B_b=\{\tilde f^r \tilde e^sb:
r,s\in \bz^+\}.$$ We say that $B_b$ is cyclic on $b$. If $\tilde
eb=0$, then $B_b$ has no branch points.

 Given branched crystals
$(B,\wt,\tilde e,\tilde f,\varepsilon,\varphi)$ and
$(B',\wt,\tilde e,\tilde f,\varepsilon,\varphi)$ a map $\Phi:B\sqcup\{0\}\to
B'\sqcup\{0\}$ is called a morphism of crystals if
\begin{enumerit}
\item[(i)] $\Phi(0)=0$,
\item[(ii)] if  $b\in B$ and $\Phi(b)\in B'$, then $\wt(\Phi(b))=\wt(b)$ and $\varepsilon(\Phi(b)) = \varepsilon(b)$,
\item[(iii)] if $b\in B$ and $\Phi(b),\Phi(\tilde fb)\in B'$ (resp. $\Phi(\tilde eb)\in B'$), then $\Phi(\tilde fb) = \tilde f\Phi(b)$ (resp. $\Phi(\tilde eb)= \tilde e\Phi(b)$).
\end{enumerit}
 A  morphism is strict if it satisfies
\begin{enumerit}
\item[(iii')]  $\Phi(\tilde fb) = \tilde f\Phi(b)$ and
 $\Phi(\tilde eb) = \tilde e\Phi(b)$ for all $b\in B$.
\end{enumerit}
An injective morphism is called an embedding and an isomorphism is a bijective morphism.

\subl{} Given $b,b'\in B$, say that $b\sim b'$ iff $B_b\cap
B_{b'}\ne\emptyset$. It follows easily  from Definition
\ref{crystal} (ii) that  if $\tilde eb=\tilde eb'=0$ then $b\sim
b'$ iff $b=b'$.
\begin{prop} Let $b,b'\in B$. Then $b\sim b'$ iff there exist
$l,l'\in \bz^+$  such that $\tilde e^lb=\tilde
e^{l'}b'\in B.$ Further,  $\sim$ is an equivalence relation on
$B$.\end{prop}
\begin{pf} The first statement is obvious from Defintion
\ref{crystal}(ii). To prove the second it is only neccesary to
check that $\sim$ is transitive.  Choose nonnegative integers
$l_i$, $1\le i\le 4$, such that
 $$\tilde e^{l_1}b_1=\tilde e^{l_2}b_2,\ \  \tilde e^{l_3}b_2=\tilde e^{l_4}b_3.$$
 If $l_2\ge l_3$, then $$\tilde
e^{l_4+l_2-l_3}b_3 =  \tilde e^{l_2}b_2 =\tilde e^{l_1}b_1 ,$$ and hence
$B_{b_1}\cap B_{b_3}\ne\emptyset$. The case $l_2<l_3$ is similar
and the lemma is proved.
\end{pf}
An equivalence class of $B$ with respect to $\sim$ is a subcrystal
and  is called an indecomposable component of $B$.
 A crystal
$B$ is indecomposable if $b\sim b'$ for all $b,b'\in B$. Thus, an
indecomposable component of $B$  is a maximal indecomposable
subcrystal of $B$.

\begin{cor}\label{uniquehw}\hfill
\begin{enumerit}
\item[(i)] An indecomposable component of a branched crystal contains at most one
element $b$ such that $\tilde eb=0$.
\item[(ii)] Suppose that $B_i$, $i\in I$, are a
family of indecomposable subcrystals of a branched crystal $B$
such that $B=\cup_{i\in I} B_i$ and
 $B_i\cap B_j=\emptyset$ if $i,j\in I$,
$i\ne j$. Then $B_i$, $i\in I$,  are the indecomposable components
of $B$.
\end{enumerit}
\end{cor}

\subl{}
\begin{prop}\label{propac} Let $B$ be a branched  crystal
and  assume that $B^{br}\ne\emptyset$. Let $b\in B^{br}$, then
\begin{enumerit}
\item[(i)]$\wt(b)\le -2$, \item[(ii)] $\tilde e^l b$ is not a
branch point for all $0< l\le -\wt (b)-1$ and  $\tilde
e^{-\wt(b)-1}b\neq 0$, \item[(iii)]$\tilde e^{-\wt(b)} b=0$,
\item[(iv)] $\tilde f^lb\neq 0 \quad\forall\ l\ge 0.$
\end{enumerit}
\end{prop}

\begin{pf}
If $b$ is a branch point then $\tilde eb\in B$ and by Definition
\ref{bc}(i),(iv),(v) we have
\begin{equation*}\label{wt}0\le \varepsilon(\tilde eb)=-\wt(\tilde
eb)= -\wt(b)-2\end{equation*} proving (i). If $\wt(b)=-2$, then
$\wt(\tilde eb)=0=\varepsilon(\tilde eb)$. The first equality
together with (i) implies that $\tilde eb\notin B^{br}$ and the
second equality implies that  $\tilde e^2b=0$ by Definition
\ref{bc}(iv). Suppose now that  $\wt(b)<-2$. It follows by induction
on $l$ that $\varepsilon(\tilde e^{l}b)=\varepsilon(\tilde
eb)-(l-1)=-\wt(b)-1-l> 0$ for $0<l<-\wt(b)-1$. Hence $\tilde e^lb\notin B^{br}$ and
$\tilde e^{l+1}b\ne 0$. Further, $\varepsilon(\tilde
e^{-\wt(b)-1}b) = 0$ and $\wt(\tilde e^{-\wt(b)-1}b) = -\wt(b) -2>
0$. This proves (ii). It also follows that $\tilde e^{-\wt(b)}
b=0$. If $l\in\bz^+$ is minimal such that $\tilde f^lb=0$ then by
Definition \ref{bc} (v) we have $l=\varphi(b)+1=\wt(b)+1$. But
this is impossible since by (i) the right hand side is negative.
\end{pf}

\begin{cor} \hfill
\begin{enumerit}
\item An indecomposable branched crystal has at most one branch point.
\item Every indecomposable crystal has one element $b$ such that $\tilde eb=0$.
\item If $\varphi(b)\ge 0$, then $\tilde e^lb\notin B^{br}$ for all $l\ge 0$.
\end{enumerit}
\end{cor}

\begin{pf} Let $B$ be an indecomposable branched crystal and
assume  that $b_1,b_2\in B^{br}$, $b_1\ne b_2$.
It follows from  Proposition \ref{propac}(iii) and Corollary
\ref{uniquehw}   that $$\tilde e^{-\wt(b_1)-1}b_1=\tilde
e^{-\wt(b_2)-1}b_2 $$and hence $\wt(b_1)=\wt(b_2)$. Applying
$\tilde f$ and using Proposition \ref{propac}(ii) we get $\tilde e
b_1=\tilde e b_2$. But this is impossible by Definition
\ref{bc}(iii). This proves (i). If $B^{br}\ne\emptyset$, (ii)
follows from Proposition\ref{propac}(iii). If $B^{br}=\emptyset$,
then  $\varepsilon(b)=\varepsilon(\tilde eb)+1 >0$ for all $b\in
B$ such that $\tilde eb\in B$.  Hence $\tilde
e^{\varepsilon(b)+1}b=0$ for all $b\in B$. To prove (iii), note
that $\varphi(b)<0$ for all $b\in B^{br}$. Therefore, it suffices
to show that if $\varphi(b)\ge 0$ then $\varphi(\tilde e^lb)\ge 0$
for all $l\ge0$ such that $\tilde e^lb\in B$. This is easily done
inductively.
\end{pf}

\subl{} The main result of this section is:
\begin{thm}\label{class} Let $B$ be an indecomposable branched crystal. Then $B$
is isomorphic to a branched crystal in one of the four infinite
families defined below. Let $r,s\in\bz, r\ge 0$.

$B(V(r))= \{b_{j}: 0\le j\le r\}$,  $$ \tilde f b_j = b_{j+1},
\qquad \tilde e b_j = b_{j-1}, \qquad \wt(b_j) = r-2j,\qquad
\varepsilon(b_j) = j,\qquad \varphi(b_j)= r-j.$$

\vskip 12pt

 $B(M(s))=\{b_j:j\ge 0\},$
\begin{gather*}
\tilde f b_j = b_{j+1}, \qquad \tilde e b_j = b_{j-1}, \qquad
\wt(b_j) = s-2j,\qquad \varepsilon(b_j) = j,\qquad \varphi(b_j)
=s-j
\end{gather*}

\vskip 24pt
 $B(T(r)) = \{b_{(j)}:j\ge 0\}\sqcup\{b_j:j\ge 0\},$
where $\tilde eb_0=b_{(r)},$ and $$\tilde f b_{(j)} = b_{(j+1)},
\qquad \tilde e b_{(j)} = b_{(j-1)}, \qquad \wt(b_{(j)}) =
r-2j,\qquad \varepsilon(b_{(j)}) =j,\qquad \varphi(b_{(j)})= r-j,
$$
 $$\tilde
f b_j = b_{j+1}, \qquad \tilde e b_{j+1} =  b_{j}, \qquad \wt(b_j)
= -r-2-2j,\qquad \varepsilon(b_j) = j,\qquad  \varphi(b_j)=
-r-2-j.$$

\vskip 24pt

 $B(M(r)^{\sigma}) = \{b_{(j)}:0\le j\le
r\}\sqcup\{b_j:j\ge 0\},$ where $\tilde eb_0=b_{(r)},$ and
$$\tilde f b_{(j)} = b_{(j+1)}, \qquad \tilde e b_{(j)} =
b_{(j-1)}, \qquad \wt(b_{(j)}) = r-2j,\qquad  \varepsilon(b_{(j)})
=j,\qquad \varphi(b_{(j)})= r-j. $$
 $$\tilde f
b_j = b_{j+1}, \qquad \tilde e b_{j+1} =  b_{j}, \qquad \wt(b_j) =
-r-2-2j,\qquad \varepsilon(b_j) = j,\qquad \varphi(b_j)= -r-2-j.
 $$
\end{thm}

 \begin{pf} It is not hard to verify that the elements of the families given in the
 theorem satisfy the conditions for  branched crystals. Now let $B$ be an indecomposable branched crystal. Suppose
 first that $B^{br}=\emptyset$ and choose $b\in B$ such that
 $\tilde eb=0$. Then, $B_{b}=\{\tilde f^r b:r\in\bz^+\}$ is a
subcrystal of $B$.  Since $\varepsilon(\tilde f^{r} b)=r$ for all
$r\in\bz^+$ such that $\tilde f^r b\in B$, it follows that $B_b$ is isomorphic to $B(M(\wt(b)))$
if $\tilde f^r b\ne 0$ for all $r\in\bz^+$, and  to
$B(V(\wt(b)))$ if $\tilde f^{\wt(b)+1}b=0$. It remains to show
that $B=B_b$. But this is clear, for given  $b'\in B$ there exists
$s\in\bz^+$ with $\tilde e^sb'=0$ or equivalently  $\tilde
e^{s-1}b'=b$ and hence $\tilde f^{s-1}b=b'$. If  $B^{br}=\{b\}$,
set $\tilde e^{-\wt(b)-1}b=\bar b$ and let $b'\in B\backslash\{b\}$.
If $\tilde e^sb'=b$ for some $s\in\bz^+$, then $b'=\tilde f^sb\in
B_b$. Otherwise $\tilde e^sb'\ne b$ for all $s\in\bz^+$ and
$\tilde e^kb'=\bar b$ for some $k\in\bz^+$. Hence $b'=\tilde
f^k \bar b\in B_b$. It is not hard to see now that $B_{b}$ is
isomorphic to $B(T(-\wt(b)-2))$ or to $B(M(-\wt(b)-2)^\sigma)$.
\end{pf}

We conclude this section with some additional results which are
needed later.

\subl{} Assume that $\{B_j\}_{j\in J}$ is a family of branched
crystals.  There exists an obvious structure of a branched crystal
on the disjoint union $\sqcup_{j\in J} B_j$ which we denote by
$\oplus_{j\in J} B_j$. Namely the maps $\wt$, $\tilde e$, $\tilde
f, \varepsilon, \varphi$ on $\oplus_{j\in J} B_j$ are defined by
requiring their restriction to $B_j$ to be the corresponding map
on $B_j$. The canonical inclusion $\iota_j: B_j\to \oplus_{j\in J}
B_j$ is a strict embedding of branched crystals.
 This proves,
\begin{prop}\label{dsdfsc} Let $B$ be a branched crystal  and let
$B_j$, $j\in J$ be the indecomposable components of $B$. Then,
$B=\oplus_{j\in J}B_j$. Furthermore, if $\Phi:B\to \oplus_{l\in L}
B'_l$ is an isomorphism of branched crystals with $B'_l$
indecomposable, then there exists a bijective map $\phi: J\to L$
such that $\Phi|_{B_j}: B_j\to B'_{\phi(j)}$, $j\in J$ is an
isomorphism of branched crystals.\end{prop}

\subl{} Given a branched crystal $B$, set  $$B^{hw}=\{b\in
B\backslash (\cup_{b'\in B^{br}}B_{b'}):\tilde eb=0\}$$ and
$$B^{br,\sigma}=\{b\in B^{br}:\tilde f\tilde eb=0\}.$$ The
following is easily seen as a consequence of Theorem \ref{class}.
\begin{prop}\label{indeccomp} Let $B$ be a branched crystal.
 The indecomposable components
of $B$ are  $B_{\bob}$, where $$\bob\in B^{hw} \cup B^{br}.$$
 Moreover we have isomorphisms of branched crystals,
\begin{align*}
&B_{\bob}\cong B(M(\wt(\bob))),\ \ {\text{if}}\ \ \bob\in B^{hw},\
\ \tilde f B_{\bob}\subset B,\\ & B_{\bob}\cong B(V(\wt(\bob))),\ \
{\text{if}}\ \ \bob\in B^{hw},\ \ \tilde f^{\wt(b)+1}\bob= 0,\\
 & B_{\bob}\cong B(T(-\wt(\bob)-2)), \ \ {\text{if}}\ \
\bob\in B^{br}\backslash B^{br,\sigma},\\
 & B_{\bob}\cong B(M(-\wt(\bob)-2)^\sigma),
\ \ {\text{if}}\ \  \bob\in B^{br,\sigma}.
\end{align*}
\end{prop}

\subl{} \begin{lem}\label{disj} Let $B$ be a branched crystal and
let $b,b'\in B$. Then $B_b\cap B_{b'}\ne\emptyset$ iff $B_{b'}\subset
B_{b}$ or vice versa.
\end{lem}

\begin{pf} Suppose  $B_{b}\cap B_{b'}\ne \emptyset$. Let $a=\tilde e^{\varepsilon(b)}b$ and $a'=\tilde e^{\varepsilon(b')}b'$. Then $B_b=B_a$ and $B_{b'}=B_{a'}$.  We have three cases:
\begin{enumerit}
\item[(i)] $a,a'\in B^{br}$,
\item[(ii)] $\tilde ea=\tilde ea'=0$,
\item[(iii)] $a\in B^{br}$  and $\tilde ea'=0$ or vice versa.
\end{enumerit}

It follows from Corollary \ref{propac}(i) (resp. Corollary \ref{uniquehw})  that $a=a'$ in case (i) (resp. case (ii)). Similarly in case (iii) we have  $\tilde e^{-\wt(a)-1}a=a'$, hence $B_{a'}\subset B_a$.
\end{pf}

\section{A tensor product rule for branched crystals}

\subl{} We now define an analogue for branched crystals of
Kashiwara's tensor product rule for crystals. Extending \cite{Ze},
we also define $\psi:B\to \mb Z$ as follows:
\begin{equation*}
\psi(b)
=\max\{\varepsilon(b),\varepsilon(b)-\varphi(b)-1\}.\end{equation*}
\nc\bF{\bold F}
\begin{defn}\label{tcd}
Given two branched crystals $B,B'$, let $B\otimes B'$ be the set
$B\times B'$ equipped with the maps $\tilde e,\tilde
f:(B\sqcup\{0\})\times (B'\sqcup\{0\})\to (B\times B')\sqcup
\{0\}$ and $\wt:B\times B'\to \bz$ defined below.
\begin{gather*}
\wt(b\otimes b') = \wt(b)+\wt(b')
\end{gather*}
\begin{align*}
\tilde f(b\otimes b')&= \tilde fb\otimes b',\text{ if } b,b'\text{
satisfy } \bF1,\\ &= b\otimes \tilde f b',\text{ if } b,b'\text{
satisfy } \bF1',\\ &= \tilde f^{\varphi(b)+1}b\otimes \tilde
e^{\varphi(b)} b',\text{ if } b,b'\text{ satisfy } \bF2,\\ \tilde
e(b\otimes b')&= \tilde eb\otimes b',\text{ if } b,b'\text{
satisfy } \be1,\\ &= b\otimes \tilde e b',\text{ if } b,b'\text{
satisfy } \be1',\\ &= \tilde e^{\varphi(b')+1}b\otimes \tilde
f^{\varphi(b')} b',\text{ if } b,b'\text{ satisfy } \be2,\\ &=
\tilde f^{-\wt(b)-\wt(b')-2}b\otimes \tilde e^{-\wt(b)-\wt(b')-1}
b',\text{ if } b,b'\text{ satisfy } \be2',\\ &=0 \text{
otherwise,}
\end{align*}
where we say that $b,b'$ satisfy
\begin{enumerit}
\item[(i)] the condition $\bF1$ if $\varphi(b)<0$ or
$\psi(b')<\varphi(b)$,
\item[(ii)] the condition $\bF1'$ if
$\tilde fb'\ne 0$ and $\psi(b')\ge \varphi(b)\ge 0$,
\item[(iii)] the  condition $\bF2$ if $
\tilde fb'= 0$ and $\psi(b')\ge \varphi(b)\ge 0$, \item[(iv)] the
condition $\be1$ if either $\wt(b)<\varphi(b)<-1$ or $\psi(b')\le
\varphi(b)$,
\item[(v)] the condition $\be1'$ if $\psi(b')> \varphi(b)\ge 0$,
\item[(vi)] the condition $\be2$ if  either $\varphi(b)=-1, \tilde f^{\varphi(b')+1}b'=0$, and
$0\le\varphi(b')\le\varepsilon(b)-1$ or $ b\in B^{br},\nexists
a'\in B'^{br}$ such that $b'=\tilde e^l a'$ for some $l\ge 0$, and
$0\le\varphi(b')\le-\wt(b)-2$,
\item[(vii)] the condition $\be2'$
if either  $\varphi(b)=-1,\tilde eb\in B$, and
$-\wt(b)-1\le \varphi(b')\le -\wt(b)-2+\varepsilon(b')$,
or $\varphi(b)<0,\varepsilon(b)=0$, and
$-\wt(b)-1\le \varphi(b')\le -\wt(b)-2+\varepsilon(b')$.
\end{enumerit}

Here we understand $0\otimes b'= b\otimes 0=0$ and define $\tilde e0=\tilde f0=0$. Finally, set

\begin{align*}
\varepsilon(b\otimes b')&=
\infty,  \text{ if } \tilde f\tilde e ^{s+1}(b\otimes b')= \tilde e^s(b\otimes b') \text{ for all } s\ge 0,\\
&= \min\{s\in \bz^+: \tilde f\tilde e ^{s+1}(b\otimes b')\ne \tilde e^s(b\otimes b')\},  \text{ otherwise }
\end{align*}
and
$$\varphi(b\otimes b')=\wt(b\otimes b')+\varepsilon(b\otimes b').$$
\end{defn}

\subl{} The next few subsections are devoted to the proof the following theorem.

\begin{thm}\label{tensor}  Let $B,B'$ be branched crystals and assume that $B_i$, $i\in
I$ and $B'_j$,  $j\in J$ are the indecomposable components of $B$ and $B'$,
 respectively. Then $B\otimes B'$ is a branched crystal and $B\otimes B' = \oplus_{i,j} B_i\otimes B'_j$.
\end{thm}

Assuming that $B\otimes B'$ is a branched crystal the second statement is proved as follows. Let $\tilde e_{B,B'}: (B\times B')\sqcup\{0\} \to (B\times B')\sqcup\{0\}$ and $\tilde e_{B_i,B'_j}: (B_i\times B'_j)\sqcup\{0\} \to (B_i\times B'_j)\sqcup\{0\}$ be defined as above. It is clear that for all $(b,b')\in B_i\times B'_j$ we have $\tilde e_{B,B'}(b\otimes b')= \tilde e_{B_i,B'_j}(b\otimes b')\in (B_i\times B'_j)\sqcup\{0\}$ (and similarly for $\tilde f$). It is now immediate that $B\otimes B' = \oplus_{i,j} B_i\otimes B'_j$.

\subl{} We proceed now with the proof of Theorem \ref{tensor}. It is clear that the condition (i) of Definition \ref{bc} is satisfied.

For (ii), let $b\in B$, $b'\in B'$ be such
that $\tilde f(b\otimes b')\in B\otimes B'$. If $b,b'$ satisfy
$\bF1$ (resp. $\bF1'$) then it is obvious that $\tilde fb,b'$
(resp. $b,\tilde fb'$) satisfy $\be1$ (resp. $\be1'$). If $b,b'$
satisfy $\bF2$, we have $$\varphi(\tilde f^{\varphi(b)+1}b) = -1,
\quad \varphi(\tilde e^{\varphi(b)}b')= \varphi(b), \quad
\varepsilon(\tilde f^{\varphi(b)+1}b) =
\varepsilon(b)+\varphi(b)+1.$$ Here we used Corollary \ref{propac}(iii). Hence $\tilde f^{\varphi(b)+1}b,
\tilde e^{\varphi(b)}b'$ satisfy $\be2$ and $\tilde e\tilde
f(b\otimes b')=b\otimes b'$. This proves (ii).

Let us now check condition (iv).
Let $b\in B, b'\in B'$ and $\bar B,\bar B'$ be the indecomposable components they belong to respectively.
It is evident that $\tilde e^l(b\otimes b')\in (\bar B\times \bar B')\sqcup\{0\}$. We claim that $\varepsilon(b\otimes b')\in \bz^+$ and, therefore,  condition (iv) of Definition \ref{bc} is  satisfied. To prove the claim it suffices to check that there exists $l\in\bz^+$ such that $\tilde e^l(b\otimes b')=0$. If that was not the case we would have $\wt(\tilde e^l(b\otimes b'))=\wt(b\otimes b')+2l$ for all $l\in\bz^+$. But this is impossible since the weight function is clearly bounded from above on $\bar B\times \bar B'$.

\subl{} To prove (iii) and (v), we begin with the following
proposition which characterizes the branch points in $B\otimes
B'$.
\begin{prop}\label{brtensor}  Suppose that $B$, $B'$ are branched
crystals. Let $b\in B$, $b\in B'$. Then \begin{enumerit}
\item[(a)] $\tilde f(b\otimes b')=0$ iff $b,b'$ satisfy $\bF2$ and
$\tilde f^{\varphi(b)+1}b=0$.
 \item[(b)] Let $b,b'$ be such that  $\tilde e(b\otimes b')\ne
0$. Then $\tilde f\tilde e(b\otimes b') \ne b\otimes b'$ iff one of the following happens:
\begin{enumerit}
\item $b'\in B'^{br}$ and $b,b'$ satisfy  $\be1'$,
\item $b\in B^{br}$ and $b,b'$ satisfy $\be2$,
\item $b,b'$ satisfy $\be2'$,
\item $\varphi(b)=-\wt(b')-2\ge 0$ and  $\varphi(b')<-1$ (or
equivalently,  $\psi(b')=\varphi(b)+1$ and
$\varphi(b')<-1$).
\end{enumerit}\end{enumerit}
\end{prop}

\begin{pf} Part (a) is obvious. For part (b),
we first  prove that $\tilde f\tilde e(b\otimes b') \ne b\otimes
b'$ if (i)--(iv) happens. This is obvious in cases (i) and (ii).
For (iii) it is also obvious unless $-\wt(b)-\wt(b')-2=0$. In this
case we need to show that $\bF1'$ does not apply to $b\otimes
\tilde eb'$. This is true since $\varphi(b)<0$. In case (iv) we
see that $\be1'$ applies to $b,b'$ and, hence, we need to show
that $\bF1'$ does not apply to $b\otimes \tilde eb'$. But this is
clear since $\varphi(b')<-1$ implies $\psi(\tilde eb')=\psi(b')-2$. The converse is proved similarly.
\end{pf}
\begin{cor} Let $b\in B$, $b'\in B'$ satisfy  $\tilde e(b\otimes b')\ne
0$ and
 $\tilde f\tilde e(b\otimes b') \ne b\otimes b'$. Then $\tilde
 f^l(b\otimes b')\in B\otimes B'$ for all $l\ge 0$.\end{cor}
 \begin{pf} Suppose that $b,b'$ satisfy condition (b)(i) of the
 proposition. Then, by Proposition \ref{propac} we see that $\tilde
 f^lb'\in B'$ for all $l\ge 0$. Since $\psi(f^lb')>\psi(b)$ it
 follows that $b$ and $f^lb'$ satisfy $\bF1'$ for all $l\ge 0$ and
 hence the corollary follows in this case.  Suppose
 now that $b,b'$ satisfy $\be2$. Then, $\tilde f^l b\in B$ for all
 $l>0$ and $\tilde f^lb$, $b'$ satisfy $\bF1$ and hence $\tilde f^l(b\otimes
 b')=\tilde f^l b\otimes b'\in B$. A similar computation
 proves the result if $b,b'$ satisfy condition (b)(iii).
Finally, if we have  $\varphi(b)=-\wt(b')-2\ge 0$ and
$\varphi(b')<-1$, then by Definition \ref{bc}(v)
we have that
$\tilde f^lb'\in B'$ for all $l\ge 0$ and since
$$\psi(b')=-\wt(b')-1\ge -\wt(b')-2=\varphi(b)\ge 0,$$ it follows
that $ b$ and $\tilde f^lb'$ always satisfy $\bF1'$ and the
corollary is proved.
 \end{pf}

\subl{}  We can now prove that condition (v) of Definition
\ref{bc} is satisfied. Let $b\otimes b'$ be such that $\tilde
e(b\otimes b')\ne 0$ and $\tilde f\tilde e(b\otimes b')\ne
b\otimes b' $. Suppose $b'\in B'^{br}$ and $b,b'$ satisfy $\be1'$.
Then $b,\tilde e^lb'$ satisfy $\be1'$ for all $0\le l\le -\varphi(b)-\wt(b')-1$ and $\bF1'$ for all $2\le l\le -\varphi(b)-\wt(b')-1$. Therefore,
 $$\tilde e^l (b\otimes b') = b\otimes \tilde e^{l}b'\ne 0,
\quad\text{for all}\quad 1\le l\le -\varphi(b)-\wt(b')-1$$
$$\tilde f\tilde e^{l+1} (b\otimes b') = \tilde e^{l}(b\otimes b')
\quad\text{for all}\quad 1\le l\le -\varphi(b)-\wt(b')-2.$$  Next,
using $\be 1$ and $\bF1$ we have  $$\tilde
e^{-\varphi(b)-\wt(b')-1+l}(b\otimes b') = \tilde e^lb\otimes
\tilde e^{-\varphi(b)-\wt(b')-1}b'\ne 0 \quad\text{for all}\quad
0\le l\le \varepsilon(b),$$ $$\tilde f\tilde
e^{-\varphi(b)-\wt(b')-1+l}(b\otimes b') = \tilde
e^{-\varphi(b)-\wt(b')-2+l}(b\otimes b') \quad\text{for all}\quad
1\le l\le \varepsilon(b),$$ and since $\tilde e^{\varepsilon(b)+1}
b=0$,  $$\tilde e^{-\varphi(b)-\wt(b')+\varepsilon(b)}(b\otimes
b')=0.$$ This proves that $$\varepsilon(\tilde e(b\otimes b'))=
-\wt(b)-\wt(b')-2,\ \ $$ i.e that $\varphi(\tilde e(b\otimes
b'))=0$. A similar computation takes care of the case $b\in
B^{br}$ and $b,b'$ satisfy $\be2$. Now suppose $b,b'$ satisfy
$\be2'$, then $$\tilde e (b\otimes b')= \tilde
f^{-\wt(b')-\wt(b)-2}b\otimes \tilde e^{-\wt(b')-\wt(b)-1}b',$$
and using $\be1$ and $\bF1$,
 it follows that
$$\tilde e^{l+1}(b\otimes b') = \tilde e^l\tilde
f^{-\wt(b')-\wt(b)-2}b\otimes \tilde e^{-\wt(b')-\wt(b)-1}b'
\quad\text{for all}\quad 0\le l\le -\wt(b')-\wt(b)-2,$$ $$\tilde
f\tilde e^{l+1}(b\otimes b') =\tilde e^{l}(b\otimes b')
\quad\text{for all}\quad 1\le l\le -\wt(b')-\wt(b)-2.$$ Since none
of the rules $\be1,\be1 ', \be2, \be2'$ apply to $b\otimes \tilde
e^{-\wt(b')-\wt(b)-1}b'$ we conclude that $\varepsilon(\tilde
e(b\otimes b'))= -\wt(b')-\wt(b)-2$, i.e.   $\varphi(\tilde
e(b\otimes b')) =0$. Finally, consider the case
$\varphi(b)=-\wt(b')-2\ge 0$ and $\varphi(b')<-1$. This time, we
find that  $\tilde e(b\otimes b')=b\otimes \tilde eb'$ and  using $\be1$ and $\bF1$ we have
$$\tilde e^{l+1}(b\otimes b') = \tilde e^lb\otimes \tilde eb'\ne 0 \quad\text{for all}\quad 0\le
l\le \varepsilon(b),$$
$$\tilde f\tilde e^{l+1}(b\otimes b') = \tilde e^{l}(b\otimes b') \quad\text{for all}\quad 1\le l\le \varepsilon(b),$$ and $\tilde e^{\varepsilon(b)+2}(b\otimes b')=0$.

Hence $\varepsilon(\tilde
e(b\otimes b')) = \varepsilon(b)$. We must check that
$\varepsilon(b) = -\wt(b)-\wt(b')-2$, but this follows from
$\varphi(b')<-1$ and $\psi(b')=\varphi(b)+1$.

To complete the proof of (v), we must show that if $l$ is minimal such that $\tilde
f^l(b\otimes b')=0$ then $l=\varphi(b\otimes b')+1$.
It clearly suffices to consider the case $l=1$. By Proposition \ref{brtensor} we know that $b,b'$ satisfy $\bF2$. Then using $\be1'$ and $\bF1'$ we see that
$$\tilde e^l(b\otimes b') = b\otimes \tilde e^lb'\ne 0, \quad\text{for all}\quad 0\le l\le \varepsilon(b')-\varphi(b), $$
and
$$\tilde f\tilde e^l(b\otimes b') = \tilde e^{l-1}(b\otimes b'), \quad\text{for all}\quad 1\le l\le \varepsilon(b')-\varphi(b).$$
Then, using $\be1$ and $\bF1$ we get
$$\tilde e^{\varepsilon(b')-\varphi(b)+l}(b\otimes b') = \tilde e^lb\otimes \tilde e^{\varepsilon(b')-\varphi(b)}b'\ne 0, \quad\text{for all}\quad 0\le l\le \varepsilon(b), $$
$$\tilde f\tilde e^{\varepsilon(b')-\varphi(b)+l}(b\otimes b') = \tilde e^{\varepsilon(b')-\varphi(b)+l-1}(b\otimes b'), \quad\text{for all}\quad 1\le l\le \varepsilon(b),$$
and $\tilde e^{\varepsilon(b')-\varphi(b)+\varepsilon(b)+1}(b\otimes b')=0$. Hence, $\varepsilon(b\otimes b')=\varepsilon(b')-\varphi(b)+\varepsilon(b) = -\wt(b')-\wt(b)$, i.e., $\varphi(b\otimes b')=0$.

\subl{} Finally, we must prove that condition (iii) of Definition \ref{bc} holds:
\begin{prop} Suppose that $b_i\in B$, $b_i'\in B'$ are distinct elements such that
$\tilde e(b_i\otimes b_i')\in B\otimes B'$ and $\tilde f\tilde
e(b_i\otimes b_i')\ne  (b_i\otimes b_i')$ for $i=1,2$. Then
$$\tilde e(b_1\otimes b_1')\ne \tilde e(b_2\otimes
b_2').$$\end{prop}
\begin{pf} Assume $\tilde e(b_1\otimes b_1')= \tilde e(b_2\otimes
b_2')$. We  consider four  cases depending on the various
possibilities for the pairs $b_1,b_1'$ given Proposition
\ref{brtensor}.

{\bf Case 1.} $b_1,b_1'$ satisfy $\be1'$ and $b_1'\in B'^{br}$.
If $b_2,b_2'$ also  satisfy $\be1'$ and $b_2'\in B'^{br}$  then
we get $b_1=b_2$  and $\tilde eb_1'=\tilde eb_2'$ which implies
that $b_1'=b_2'$. If $b_2,b_2'$ satisfy $\be2$  with $b_2\in B^{br}$, then we get
$$b_1=\tilde e^{\varphi(b_2')+1}b_2,\ \ \tilde eb_1'=\tilde
f^{\varphi(b_2')}b_2',$$ i.e., $b_2'= \tilde
e^{\varphi(b_2')+1}b_1'$. But
 $\varphi(b_2)=\wt(b_2)\le -2$ by Proposition \ref{propac},
 so we get a contradiction to the fact that $b_2,b_2'$ satisfy
 $\be2$. Next, suppose that $b_2,b_2'$ satisfy $\be2'$. This
 gives,
 $$b_1=\tilde f^{-(\wt(b_2)+\wt(b_2')+2)}b_2,$$ which gives
 $\varphi(b_1)= \varphi(b_2)+\wt(b_2)+\wt(b_2')+2<0$, contradicting
 the fact that $\varphi(b_1)\ge 0$.
 Finally, suppose that $b_2,b_2'$ satisfy  $\varphi(b_2)=-\wt(b_2)-2\ge 0$ and
  $\varphi(b_2')<-1$. In particular, this means that $b_2,b_2'$
  satisfy $\be1'$ and hence we get $b_1=b_2$ and $\tilde
  eb_1'=\tilde eb_2'$. Since $b_1'\in B'^{br}$ this implies that
  $b_2'\notin B'^{br}$ and so $\tilde f\tilde eb_1'=b_2'$, which
  gives $$\varphi(b_2')=\wt(b_1')+\varepsilon(\tilde f\tilde e
  b_1')=\wt(b_1')+\varepsilon(\tilde eb_1')+1 =-1,$$ since
  $\varepsilon(\tilde eb_1')=-\wt(b_1')-2$ which contradicts
  $\varphi(b_2')<-1$.

  {\bf Case 2.} $b_1$, $b_1'$ satisfy $\be2$ and $b_1\in B^{br}$.
  If $b_2$ and $b_2'$ also satisfy $\be2$ and $b_2\in B^{br}$,
  assume without loss of generality that $\varphi(b_1')\ge
  \varphi(b_2')$. Since,
  $\tilde e^{\varphi(b_1')+1}b_1=\tilde e^{\varphi(b_2')+1}b_2$,
  we get by using Proposition \ref{propac} that $\tilde e b_1=\tilde
  f^{\varphi(b_1')-\varphi(b_2')}\tilde eb_2$. But this means
  that,
  $$0=\varphi(\tilde   eb_1)=-\varphi(b_1')+\varphi(b_2'),$$
   which implies $b_1=b_2$. Since we also have $\tilde f^{\varphi(b_1')-\varphi(b_2')}b_1'=b_2'$, it follows that
   $b_1'=b_2'$.

  Assume next that $b_2$ and $b_2'$ satisfy $\be2'$.
  We get $\tilde e^{\varphi(b_1')+1}b_1=\tilde
  f^{-\wt(b_2')-\wt(b_2)-2}b_2$ and hence $b_2=\tilde
  e^{\varphi(b_1')-\wt(b_2')-\wt(b_2)-1}b_1.$ This implies that
  $$\varphi(b_2)=\varphi(\tilde
  eb_1)+\varphi(b_1')-\wt(b_2')-\wt(b_2)-2\ge 0, $$ which is again a contradiction.
  Finally, let us assume that $\varphi(b_2)=-\wt(b_2')-2\ge 0$ and
  $\varphi(b_2')<-1$, then $\tilde e(b_2\otimes b_2')=b_2\otimes
  \tilde eb_2'$ and so we get $\tilde
  e^{\varphi(b_1')+1}b_1=b_2$, which gives
  $\varphi(b_2)=\varphi(b_1')<0$ which is again a contradiction.

  {\bf Case 3.} $b_1$ and $b_1'$ satisfy $\be2'$.  If $b_2$ and
  $b_2'$ also satisfy $\be2'$, then we get $\tilde
  f^{-\wt(b_1)-\wt(b_1')-2}b_1= \tilde f^{-\wt(b_2)-\wt(b_2')-2}b_2$.
  Since $\wt(b_1)+\wt(b_1')=\wt(b_2)+\wt(b_2') (=l)$, it follows
  immediately that $b_1=b_2$. Since
  $\tilde e^{l-1}b_1'=\tilde e^{l-1}b_2'$ and $\varphi(b_i')\ge 0$, it follows from Corollary \ref{propac}(iii) that
  $b_1'=b_2'$.

If $b_2$ and $b_2'$
satisfy  $\varphi(b_2)=-\wt(b_2')-2\ge 0$ and
  $\varphi(b_2')<-1$, we get $b_2=\tilde
  f^{-\wt(b_1)-\wt(b_1')-2}b_1$. This means that $\tilde e^{-\wt(b_1)-\wt(b_1')-2}b_2=b_1$
   and that  $\tilde
  e^lb_2\notin B^{br}$ for all $0\le l\le -\wt(b_1)-\wt(b_1')-2$  and hence $\varphi(b_1)=\varphi(b_2)-\wt(b_1)-\wt(b_1')-2\ge 0$
  which is a contradiction.

  {\bf Case 4} The remaining case is when $b_i$ and $b_i'$ satisfy
  $\varphi(b_i)=-\wt(b_i')-2\ge 0$ and
  $\varphi(b_i')<-1$, $i=1,2$. Again we get $b_1=b_2$ and $\tilde
  eb_1'=\tilde eb_2'$. This means that  the only time that it is not obvious that $b_1'=b_2'$ is when
   $b_1'\in B'^{br}$ and  $b_2'\notin
  B'^{br}$(or vice versa). But then we get $b_2'=\tilde f\tilde
  eb_1'$ and so $\varphi(b_2')=-1$ since $\varphi(\tilde
  eb_1')=0$. But this is again a contradiction.
\end{pf}

The proof of Theorem \ref{tensor} is now complete.

\subl{}
\begin{prop}\label{hwtensor}
Let $B,B'$ be branched crystals and let $b\in B$, $b'\in B'$.
Then
 $\tilde e(b\otimes b')= 0$ iff one of the following holds:
\begin{enumerit}
\item $\varepsilon(b)=0, b\notin B^{br}$ (resp. $\varepsilon(b')=0, b'\notin B'^{br}$), and $b,b'$ satisfy $\be1$ (resp.
$\be1'$),or  \item $b,b'$ do not satisfy  $\be1$, $\be2$, $\be1'$,
$\be2'$.
\end{enumerit}
Moreover, if (ii) holds, then either $\varepsilon(b)=0$ or
$\varphi(b)=-1$\end{prop}
\begin{pf} Observe that if $b,b'$ satisfy $\be1$ (resp. $\be1'$) then $\tilde
e(b\otimes b')=0$ iff $\tilde eb=0$ (resp. $\tilde eb'=0$).
It is also easy to see that if $b,b'$ satisfy none of the conditions
$\be1,\be1',\be2,\be2'$, then either $\varepsilon(b)=0$ or
$\varphi(b)=-1$. It remains to show that if $b,b'$ satisfy $\be2$
or $\be2'$, then $\tilde e(b\otimes b')\ne 0$. If $b,b'$ satisfy
$\be2$, observe that the upper bound for $\varphi(b')$ implies that
$\tilde e^{\varphi(b')+1}b\ne 0$ and we are done. If $b,b'$
satisfy $\be2'$ then, since $\varphi(b)<0$, we have $\tilde f^lb\ne
0$ for all $l\ge 0$. We need to show that $\tilde
e^{-\wt(b)-\wt(b')-1}b'\ne 0$. It suffices to check that
$-\wt(b)-\wt(b')-1\le \varepsilon(b')$. But this is equivalent to
$\varphi(b')\ge -\wt(b)-1$.
\end{pf}

\subl{}\label{tpdc}

Let $B,B'$ be indecomposable branched crystals and $\bb = B\otimes
B'$. We conclude this section by writing down  the indecomposable
components of $\bb$ in a case by case fashion. By Proposition
\ref{indeccomp}  it suffices to compute the sets
$\bb^{br},\bb^{br,\sigma},\bb^{hw}$  which is done by using
Propositions \ref{brtensor} and \ref{hwtensor} and Lemma
\ref{disj}. In the following we use the notation of Theorem
\ref{class}. Also we understand $p\in\bz^+$ and use the convention
$B(M(-1))=B(T(-1))$ when convenient.  \vskip 12pt

\noindent {\bf Case 1.} Suppose that $B=B(M(s)), B'=B(V(r))$.

\begin{align*}
& \bb_{b_{s+1}\otimes b'_{r-p}}\cong B(T(r+s-2p)), 2s+2\le 2p\le r+s+1, \text{ if } s\ge 0,\\
& \bb_{b_0\otimes b'_{r+s+1-p}}\cong B(T(r+s-2p)), 0\le 2p\le r+s+1, \text{ if } s< 0,\\
& \bb_{b_0\otimes b'_{p}}\cong B(M(r+s-2p)), 0\le p\le \min\{r,s\},\\
& \bb_{b_0\otimes b'_{p}}\cong B(M(r+s-2p)), r+s+2\le p\le r.
\end{align*}

\noindent {\bf Case 2.} $B=B(T(s)), B'=B(V(r))$.

\begin{align*}
& \bb_{b_{(s+1)}\otimes b'_{r-p}}\cong B(T(r+s-2p)), 2s+2\le 2p\le r+s+1,\\
& \bb_{b_0\otimes b'_{r-p}}\cong B(T(r+s-2p)), 0\le 2p\le r+s+1, \text{ if } s<r,\\
& \bb_{b_0\otimes b'_{r-p}}\cong B(T(r+s-2p)), 0\le p\le r, \text{ if } s\ge r.
\end{align*}

\noindent {\bf Case 3.} $B=B(M(s)^{\sigma}), B'=B(V(r))$.

\begin{align*}
& \bb_{b_0\otimes b'_{r-p}}\cong B(T(r+s-2p)), 2s+2\le 2p\le r+s+1,\\
& \bb_{b_0\otimes b'_{r-p}}\cong B(M(r+s-2p)^{\sigma}), 0\le p\le \min\{r,s\}.
\end{align*}

\noindent {\bf Case 4.} $B=B(V(r)), B'=B(M(s)^{\sigma})$.

\begin{align*}
& \bb_{b_{r+s-2p}\otimes b'_{p-s}}\cong B(T(r+s-2p)), 2s+2\le 2p\le r+s,\\
& \bb_{b_{r-p}\otimes b'_{0}}\cong B(M(r+s-2p)^{\sigma}), 0\le p\le \min\{r,s\},\\
& \bb_{b_0\otimes b'_{(r-s-1)/2}} \cong B(M(-1)).
\end{align*}

\noindent {\bf Case 5.} $B'=B(M(r)), B=B(M(s))$.

\begin{align*}
& \bb_{b_{r+s-2p}\otimes b'_{p+1}}\cong B(T(r+s-2p)), 2s+2\le 2p\le r+s,\\
& \bb_{b_{r+1}\otimes b'_{s-p}}\cong B(T(r+s-2p)), 2r+2\le 2p\le r+s+1, \text{ if } r\ge 0,\\
& \bb_{b_{0}\otimes b'_{r+s+1-p}}\cong B(T(r+s-2p)), 0\le 2p\le r+s+1, \text{ if } r< 0,\\
& \bb_{b_0\otimes b'_p}\cong B(M(r+s-2p)), 0\le p\le \min\{r,s\},\\
& \bb_{b_0\otimes b'_p}\cong B(M(r+s-2p)), p\ge r+s+2, \text{ if } r< 0,\\
& \bb_{b_{r+1}\otimes b'_{p-r-1}}\cong B(M(r+s-2p)), p\ge \max\{r,s\}+1, \text{ if } r\ge 0,\\
& \bb_{b_p\otimes b'_0}\cong B(M(r+s-2p)), r+s+2\le p\le r, \text{ if } s<0,\\
& \bb_{b_0\otimes b'_{(r+s+1)/2}} \cong B(M(-1)), \text{ if } r\ge s.\\
\end{align*}

\noindent {\bf Case 6.} $B=B(M(r)), B'=B(M(s)^{\sigma})$.

\begin{align*}
& \bb_{b_{r-p}\otimes b'_{0}}\cong B(T(r+s-2p)), 0\le p\le s,\\
& \bb_{b_{r+s-2p}\otimes b'_{p-s}}\cong B(T(r+s-2p)), 2s+2\le 2p\le r+s,\\
& \bb_{b_{r+1}\otimes b'_{(s-p)}}\cong B(T(r+s-2p)), 2r+2\le 2p\le r+s+1, \text{ if } r\ge 0,\\
& \bb_{b_{0}\otimes b'_{(r+s+1-p)}}\cong B(T(r+s-2p)), 0\le 2p\le r+s+1, \text{ if } r< 0,\\
& \bb_{b_{r+1}\otimes b'_{p-r-s-2}}\cong B(M(r+s-2p)), p\ge r+s+2, \text{ if } r\ge 0,\\
& \bb_{b_{0}\otimes b'_{(p)}}\cong B(M(r+s-2p)), r+s+2\le p\le s, \text{ if } r< 0,\\
& \bb_{b_{0}\otimes b'_{p-s-1}}\cong B(M(r+s-2p)), p\ge s+1, \text{ if } r< 0,\\
& \bb_{b_0\otimes b'_{(r-s-1)/2}}\cong M(-1), \text{ if } r\ge 0.
\end{align*}

\noindent {\bf Case 7.} $B=B(M(s)^{\sigma}), B'=B(M(r))$.

\begin{align*}
& \bb_{b_0\otimes b'_{r-p}}\cong B(T(r+s-2p)), 0\le 2p\le r+s+1, \text{ if } s<r,\\
& \bb_{b_0\otimes b'_{r-p}}\cong B(T(r+s-2p)), 0\le p\le r, \text{ if } s\ge r,\\
& \bb_{b_{(r+s-2p)}\otimes b'_{p+1}}\cong B(T(r+s-2p)), 2r+2\le 2p\le r+s,\\
& \bb_{b_{0}\otimes b'_{p-s-1}}\cong B(M(r+s-2p)), p\ge r+s+2, \text{ if } r\ge 0,\\
& \bb_{b_{0}\otimes b'_{p-s-1}}\cong B(M(r+s-2p)), p\ge s+1, \text{ if } r< 0, \\
& \bb_{b_{(p)}\otimes b'_{0}}\cong B(M(r+s-2p)), s+r+2\le p\le s,\\
& \bb_{b_{(0)}\otimes b'_{(r+s+1)/2}}\cong B(M(-1)), \text{ if }  s \ge r.
\end{align*}

\noindent {\bf Case 8.} $B=B(M(r)^{\sigma}), B'=B(M(s)^{\sigma})$.

\begin{align*}
& \bb_{b_{(r+s-2p)}\otimes b'_{p-s}}\cong B(T(r+s-2p)),2s+2\le 2p\le r+s,\\
& \bb_{b_{0}\otimes b'_{(s-p)}}\cong B(T(r+s-2p)),2r+2\le 2p\le r+s+1,\\
& \bb_{b_{(r-p)}\otimes b'_{0}}\cong B(M(r+s-2p)^{\sigma}),0\le p\le \min\{r,s\},\\
& \bb_{b_{0}\otimes b'_{(p-r-1)}}\cong B(M(r+s-2p)),\max\{r,s\}+1\le p\le r+s+1,\\
& \bb_{b_{0}\otimes b'_{p-r-s-2}}\cong B(M(r+s-2p)),p\ge r+s+2,\\
& \bb_{b_{(0)}\otimes b'_{(r-s-1)/2}}\cong B(M(-1)).
\end{align*}

\noindent {\bf Case 9.} $B=B(M(s)), B'=B(T(r))$.

\begin{align*}
& \bb_{b_{s-p}\otimes b'_{0}}\cong B(T(r+s-2p)), 0\le p\le \min\{r,s\},\\
& \bb_{b_{r+s-2p}\otimes b'_{p-r}}\cong \bb_{b_{r+s-2p}\otimes b'_{(p+1)}}\cong B(T(r+s-2p)), 2r+2\le 2p\le r+s,\\
& \bb_{b_{s+1}\otimes b'_{(r-p)}}\cong B(T(r+s-2p)), 2s+2\le 2p\le r+s+1,\\
& \bb_{b_{0}\otimes b'_{(r+s+1-p)}}\cong B(T(r+s-2p)), 0\le 2p\le r+s+1, \text{ if } s<0,\\
& \bb_{b_{s+1}\otimes b'_{p-r-s-2}}\cong B(M(r+s-2p)), p\ge r+s+2, \text{ if } s\ge 0,\\
& \bb_{b_{s+1}\otimes b'_{(p-s-1)}}\cong B(M(r+s-2p)), p\ge \max\{r,s\}+1, \text{ if } s\ge 0,\\
& \bb_{b_{0}\otimes b'_{(p)}}\cong B(M(r+s-2p)), p\ge r+s+2, \text{ if } s< 0,\\
& \bb_{b_{0}\otimes b'_{(p-s-1)}}\cong B(M(r+s-2p)), p\ge s+1, \text{ if } s< 0,\\
& \bb_{b_0\otimes b'_{(s-r-1)/2}}\cong \bb_{b_0\otimes b'_{((r+s+1)/2)}}\cong B(M(-1)), \text{ if } s\ge r.
\end{align*}

\noindent {\bf Case 10.} $B=B(T(r)), B'=B(M(s))$.

\begin{align*}
& \bb_{b_0\otimes b'_{s-p}}\cong B(T(r+s-2p)), 0\le 2p\le r+s+1, \text{ if } s>r,\\
& \bb_{b_0\otimes b'_{s-p}}\cong B(T(r+s-2p)), 0\le p\le s,  \text{ if } s\le r,\\
& \bb_{b_{(r+1)}\otimes b'_{s-p}}\cong B(T(r+s-2p)), 2r+2\le 2p\le r+s+1,\\
& \bb_{b_{(r+s-2p)}\otimes b'_{p+1}}\cong B(T(r+s-2p)), 2s+2\le 2p\le r+s,\\
& \bb_{b_{0}\otimes b'_{p-r-1}}\cong B(M(r+s-2p)), p\ge r+s+2, \text{ if } s\ge 0,\\
& \bb_{b_{0}\otimes b'_{p-r-1}}\cong B(M(r+s-2p)), p\ge r+1, \text{ if } s<0,\\
& \bb_{b_{(r+1)}\otimes b'_{p-r-1}}\cong B(M(r+s-2p)), p\ge \max\{r,s\}+1,\\
& \bb_{b_{(p)}\otimes b'_{0}}\cong B(M(r+s-2p)), s+r+2\le p\le r,\\
& \bb_{b_{(0)}\otimes b'_{(r+s+1)/2}}\cong B(M(-1)), \text{ if } s \le r.
\end{align*}

\noindent {\bf Case 11.} $B=B(T(r)), B'=B(T(s))$.

\begin{align*}
& \bb_{b_{(r+s-2p)}\otimes b'_{(p+1)}}\cong \bb_{b_{(r+s-2p)}\otimes b'_{p-s}}\cong B(T(r+s-2p)), 2s+2\le 2p\le r+s,\\
& \bb_{b_{(r-p)}\otimes b'_{0}}\cong B(T(r+s-2p)), 0\le p\le s,\\
& \bb_{b_{(r+1)}\otimes b'_{(s-p)}}\cong \bb_{b_0\otimes b'_{(s-p)}}\cong B(T(r+s-2p)), 2r+2\le 2p\le r+s+1,\\
& \bb_{b_{(r+1)}\otimes b'_{(p-r-1)}}\cong \bb_{b_{0}\otimes b'_{(p-r-1)}}\cong B(M(r+s-2p)), p\ge s+1,\\
& \bb_{b_{(r+1)}\otimes b'_{p-r-1}}\cong \bb_{b_{0}\otimes b'_{p-r-1}}\cong (M(r+s-2p)), p\ge r+1,\\
& \bb_{b_{(0)}\otimes b'_{((r+s+1)/2)}} \cong \bb_{b_{(0)}\otimes b'_{(r-s-1)/2}} \cong B(M(-1)), \text{ if } r\ge s.\\
\end{align*}

\noindent {\bf Case 12.} $B=B(T(r)), B'=B(M(s)^{\sigma})$.

\begin{align*}
& \bb_{b_{(r+s-2p)}\otimes b'_{p-s}}\cong B(T(r+s-2p)), 2s+2\le 2p\le r+s,\\
& \bb_{b_{(r-p)}\otimes b'_{0}}\cong B(T(r+s-2p)), 0\le p\le \max\{r,s\},\\
& \bb_{b_{0}\otimes b'_{(s-p)}}\cong \bb_{b_{(r+1)}\otimes b'_{(s-p)}}\cong B(T(r+s-2p)), 2r+2\le 2p\le r+s+1,\\
& \bb_{b_{0}\otimes b'_{(p-r-1)}}\cong B(M(r+s-2p)), p\ge r+1, \text{ if } r\ge s,\\
& \bb_{b_{0}\otimes b'_{(p-r-1)}}\cong B(M(r+s-2p)), s+1\le p\le r+s+1, \text{ if } s> r,\\
& \bb_{b_{0}\otimes b'_{p-r-s-2}}\cong \bb_{b_{(r+1)}\otimes b'_{p-r-s-2}}\cong B(M(r+s-2p)), p\ge r+s+2,\\
& \bb_{b_{(0)}\otimes b'_{(r-s-1)/2}}\cong B(M(-1)), \text{ if } r\ge s.\\
\end{align*}

\noindent {\bf Case 13.} $B=B(M(s)^{\sigma}), B'=B(T(r))$.

\begin{align*}
& \bb_{b_{(r+s-2p)}\otimes b'_{(p+1)}}\cong \bb_{b_{(r+s-2p)}\otimes b'_{p-r}}\cong B(T(r+s-2p)), 2r+2\le 2p\le r+s,\\
& \bb_{b_{(s-p)}\otimes b'_{0}}\cong B(T(r+s-2p)), 0\le p\le r,\\
& \bb_{b_0\otimes b'_{(r-p)}}\cong B(T(r+s-2p)), 2s+2\le 2p\le r+s+1,\\
& \bb_{b_{0}\otimes b'_{(p-s-1)}}\cong B(M(r+s-2p)), p\ge r+1,\\
& \bb_{b_{0}\otimes b'_{p-s-1}}\cong (M(r+s-2p)), p\ge s+1,\\
& \bb_{b_{(0)}\otimes b'_{((r+s+1)/2)}} \cong B(M(-1)), \text{ if } r\le s.\\
\end{align*}

The cases $B=B(V(r))$ and  $B'\in\{B(V(s)), B(M(s)), B(T(s))\}$ were done in \cite{K}, \cite{Ze}. We give it here for completeness.

\noindent {\bf Case 14.} $B=B(V(r))$, $B'=B(V(s))$.

$$\bb_{b_0\otimes b'_p} \cong B(V(r+s-2p)), 0\le p\le \min\{r,s\}.$$

\noindent {\bf Case 15.} $B=B(V(r))$, $B'=B(M(s))$.

\begin{align*}
& \bb_{b_{r+s-2p}\otimes b'_{p+1}}\cong B(T(r+s-2p)), 2s+2\le 2p\le r+s,\\
& \bb_{b_0\otimes b'_p}\cong B(M(r+s-2p)), 0\le p\le \min\{r,s\}, \\
& \bb_{b_p\otimes b'_0}\cong B(M(r+s-2p)), r+s+2\le p\le r,\\
& \bb_{b_0\otimes b'_{(r+s+1)/2}} \cong B(M(-1)), \text{ if } r\ge s.\\
\end{align*}

\noindent {\bf Case 16.} $B=B(V(r))$, $B'=B(T(s))$.

\begin{align*}
& \bb_{b_{r+s-2p}\otimes b'_{(p+1)}}\cong \bb_{b_{r+s-2p}\otimes b'_{p-s}}\cong B(T(r+s-2p)), 2s+2\le 2p\le r+s,\\
& \bb_{b_{r-p}\otimes b'_{0}}\cong B(T(r+s-2p)), 0\le p\le \min\{r,s\},\\
& \bb_{b_0\otimes b'_{((r+s+1)/2)}} \cong \bb_{b_0\otimes b'_{(r-s-1)/2}} \cong B(M(-1)), \text{ if } r\ge s.\\
\end{align*}

\section{ The category $\cal O$ and branched crystals.} \label{s:repth}

We recall the definition and properties of the category $\cal O$
of modules for $\bu_q(sl_2)$ in 3.1--3.4. The relevant results can
be found in \cite {BGG}, \cite{DGK}, \cite {J2}, and \cite{So}.

\subl{}  Let $\mb C(q)$ be the field of rational functions in an
indeterminate $q$ and let $\bu_q(\lie{sl_2})$ be the quantized
enveloping algebra of $\lie{sl_2}$  over $\mb C(q)$, i.e.,  the
algebra generated by elements $e, f, k^{\pm 1}$ and relations $$
kk^{-1}=1=k^{-1}k,\ \ kek^{-1}=q^2e,\ \ kfk^{-1}=q^{-2}f,$$
$$ef-fe=\frac{k-k^{-1}}{q-q^{-1}}.$$

 It is well--known that $\bu_q(\lie{sl_2})$ is a Hopf algebra.
Let $\sigma:\bu_q(\lie{sl_2})\to \bu_q(\lie{sl_2})$ be the involutive
antiautomorphism obtained by extending the assignment
$$\sigma(e)=f,\ \ \sigma(f)=e,\ \ \sigma(k)=k^{-1},\ \
\sigma(q)=q^{-1}.$$ Let $\Omega\in\bu_q(\lie sl_2)$ be the quantum
Casimir element.

\subl{} Recall that a module $M$ of $\bu_q(\lie{sl_2})$ is said to
be a weight module of type 1 if we can write $$M=\oplus_{r\in\mb
Z} M_r,\ \ \ \  \ M_r=\{m\in M: km=q^{r}m\}.$$ Let $\cal O$ be the
category of type 1 $\bu_q(\lie {sl_2})$--modules $M$ such that
 $\text{dim} M_r<\infty$ for all $r\in\bz$,
and  such that there exists an integer $n$ depending on $M$ such
that $M_r=0$ for all $r\ge n$. It is obvious that  $\cal O$ is an
abelian category and that it is closed under taking tensor
products.
 Given any $M\in\cal{O}$ let $M^\sigma$ be the subspace of $M^*={\rm
Hom}_{\mb C(q)}(M,\mb C(q))$ consisting of elements $m^*$ such
that $m^*(M_r)=0$ for all but finitely many $r\in\bz$. The formula
$$(um^*)(m)=m^*(\sigma(u)m),\ \ u\in U_q(\lie{sl_2}),\ m^*\in
M^*,\ m\in M,$$ defines an action of $\bu_q(\lie{sl_2})$ on
$M^\sigma$ and it is easy to check that  $M^\sigma\in\cal O$. The
module $M^\sigma$ is called the restricted dual of $M$ and we say
$M$ is self--dual if $M\cong M^\sigma$. The assignment $M\mapsto
M^\sigma$ is an exact contravariant functor on $\cal{O}$ and , if
$M,N\in\cal{O}$, then $(M\otimes N)^\sigma\cong N^\sigma\otimes
M^\sigma$. Given $M\in\cal{O}$, and $c\in\mb C(q)$, let
$$M[c]=\{m\in M:(\Omega -c)^nm =0,\ \ \forall\ n>>0\}.$$ Clearly
$M[c]$ is a $\bu_q(\lie{sl_2})$--submodule of $M$, $
M=\oplus_{c\in\mb C(q)} M[c]$, and
 $M[c]\ne 0$ implies
$c= (q^{r+1}+q^{-r-1})(q-q^{-1})^{-2}$  for some $r\in\mb Z.$

\subl{} Given $r\in\mb Z$, $r\ge -1$, let $\cal{O}_r$ be the
abelian subcategory consisting of modules $M\in\cal{O}$ such that
$$M=M[(q^{r+1}+q^{-r-1})(q-q^{-1})^{-2}].$$ Since
$\sigma(\Omega)=\Omega$ it follows that  $M\in\cal{O}_r$ iff
$M^\sigma\in\cal{O}_r.$ Further, $$\cal{O}=\oplus_{r\ge
-1}\cal{O}_r$$ and  any object in $\cal{O}_r$  is isomorphic to a
direct sum of indecomposable objects in $\cal{O}_r$. The last
statement is easily seen by  an induction on
$\dim(M_r)+\dim(M_{-r-2})$, $r\ge -1$.

 \subl{} We now describe the indecomposable objects in
$\cal{O}_r$. Given $s\in\bz$, let $M(s)$ be the Verma module with
highest weight $s$ and highest weight vector $m_s$. Namely, $M(s)$
is the $\bu_q(\lie{sl_2})$--module generated by an element $m_s$
with relations: $$em_s=0,\ \ km_s=q^{s}m_s.$$ Let $V(s)$ be the
unique irreducible quotient of $M(s)$. If $r\ge 0$, then
\begin{equation*}\label{irred}
M(-r-1)\cong V(-r-1) \quad\text{ and }\quad {\dim} V(r)=r+1.
\end{equation*}
For  $r\in\bz^+$, let $T(r)$ be the $\bu_q(\lie{sl_2})$--module
generated by an element $t_r$ satisfying the relations: $$e^{r+2}
t_r=0, \ \ k t_r=q^{-r-2}t_r,\
(\Omega-(q^{r+1}+q^{-r-1})(q-q^{-1})^{-2})^2 t_r=0. $$ Note that
$\dim(T(r))_{-r-2p}=2$ if $p> 0$ and $\dim T(r)_{r-2p}=1$ if $0\le
p\le r$.

For $r\in\bz^+$, the modules $M(-r-1)$, $T(r)$, and $V(r)$ are
self dual. Let $m_r^*\in M(r)^\sigma$ be such that $m_r^*(m_r)=1$
and $m_r^*(M(r)_s)=0$ if $s\ne r$. Then
$$U_q(\lie{sl_2})m_r^*\cong V(r).$$

The following proposition is proved in \cite{So}.
\begin{prop} Let $r\in\bz^+$. The modules $M(r)$, $M(r)^\sigma$,
$M(-r-2)$, $T(r)$, and $V(r)$ are precisely the indecomposable
modules in $\cal{O}_r$ while the only indecomposable module in
$\cal{O}_{-1}$  is $M(-1)$.
\end{prop}

\begin{cor} Let $r\in\bz^+$ and $V\in\cal{O}_r$. There exist unique integers
$n_V(T(r))$, $n_V(M(r))$, etc, such that,
\begin{equation*}
 V\cong T(r)^{\oplus n_V(T(r))}\oplus M(r)^{\oplus n_V(M(r))}\oplus
 M(-r-2)^{\oplus n_V(M(-r-2))} \oplus  V(r)^{\oplus n_V(V(r))}\oplus  (M(r)^\sigma)^{\oplus n_V((M(r))^\sigma)}.
\end{equation*}
\end{cor}

\subl{}\label{tpdr} It follows that the tensor product of any two objects of
$\cal{O}$ can be written uniquely as a direct sum  of
indecomposable modules. We give these decompositions explicitly. It is enough to write down the formulas
for modules $V$ in the set $$\{V(r)\otimes V(s), \ V(r)\otimes
M(s),\ V(r)\otimes T(s),\ T(r)\otimes M(s),\ T(r)\otimes T(s),\
M(r)\otimes M^\sigma(s),\ M(r)\otimes M(s)\}$$ since the other
decompositions can be obtained by using $\sigma$. Our proof
actually establishes that the tensor product in $\cal O$ is
commutative, although this can be established as usual by using
the $R$--matrix. \vskip 6pt
 \noindent {\bf Case 1.}  If
$V=V(r)\otimes V(s)$, $r,s\ge 0$ then  it is well--known that
$V(r)\otimes V(s)\cong\bigoplus_{p=0}^{\min\{r,s\}}V(r+s-2p).$
 \noindent{\bf Case 2.}  In all other cases,  the
element $f$ acts freely on $V$ and hence we can write,
\begin{equation}\label{decomp1}V\cong\mathbf M_{<-1}\oplus\mathbf T\oplus\mathbf
M_{\ge 0},\end{equation}
 where \begin{equation}\label{decomp2}\mathbf M_{<-1}=
 \bigoplus_{s\ge 0} M(-s-2)^{\oplus n_V(M(-s-2))},\ \
 \mathbf T=\bigoplus_{r\ge -1} T(r)^{\oplus n_V(T(r))},\ \ \mathbf
 M_{\ge 0}=\bigoplus_{r\ge 0}M(r)^{\oplus
 n_V(M(r))},\end{equation}
 where we set $T(-1)=M(-1)$.
\vskip 12pt \vskip 6pt

\noindent{\bf Case 2(a)} Suppose that $V=M(r)\otimes M(s)$ and
$r,s\ge 0$. Then, $$\mathbf M_{\ge 0}=\bigoplus_{p=0}^{\min\{r,s\}}
M(r+s-2p).$$ To see this, observe that \eqref{decomp1} implies
that $$V^\sigma\cong \bigoplus_{\ell\ge 0}
(M(\ell)^\sigma)^{\oplus {n_V(M(\ell))}} \oplus\mathbf
T\oplus\mathbf M_{<-1}.$$   It is clear that the submodule
$$V^\sigma_{tor}=\{v\in V^\sigma: f^\ell v=0,\ \ \ell>0\}$$ of $V$
is isomorphic to the  submodule $\bigoplus_{\ell\ge 0}
V(\ell)^{\oplus {n_V(M(\ell))}}$ of $\bigoplus_{\ell\ge 0}
(M(\ell)^\sigma)^{\oplus {n_V(M(\ell))}}$. On the other hand, it
is not hard to see that
$$V^\sigma_{tor}\cong V(r)\otimes V(s).$$ It follows immediately that
$$n_V(M(r+s-2p))=1,\ \ 0\le p\le\min\{r,s\},$$ and $n_V(M(\ell))=0$
otherwise. Since $\dim(M(r)\otimes M(s))_{r+s-2p}=p+1$ for all
$p\ge 0$, a simple dimension counting now implies that $$\mathbf T
= \bigoplus_{p=\min\{r,s\}+1}^{[\frac{r+s+1}{2}]} T(r+s-2p),
\quad\text{ and }\quad \mathbf M_{<-1} =
\bigoplus_{p=\max\{r,s\}+1}^{\infty} M(r+s-2p).$$

\noindent{\bf Case 2(b)} If $V=V(r)\otimes M(s)$ with $r,s\ge 0$,
then
\begin{alignat*}{3}
&\mathbf M_{<-1} & = & 0,\\ &\mathbf M_{\ge 0}&  = &
\bigoplus_{p=0}^{\min\{r,s\}} M(r+s-2p),\\ &\mathbf T  &=
&\bigoplus_{p=s+1}^{[\frac{r+s+1}{2}]} T(r+s-2p), \ \text{ if } \
r>s,\\ && = &0, \ \text{ if } r\le s.
\end{alignat*}
This can be proved similarly, see also \cite{E}.

\noindent {\bf Case 2(c).} If $V\notin\{ M(r)\otimes M( s),
V(r)\otimes M(s): r,s\in\bz^+\}$, then
   $f$ acts freely on both $V$ and $V^\sigma$ and hence $\bm_{\ge
   0}=0$. This means that $V_p\cong \mathbf T_p$  for all $p\ge
-1$, and so by comparing the dimension of the weight spaces on
both sides of \eqref{decomp1} we get
 $$ n_V(T(p))=\dim V_p-\dim V_{p+2},\ \ \forall\  p\ge -1,$$
and, for $\ell \ge 0$,
  \begin{align*}
n_V(M(-\ell-2))&= \dim V_{-\ell-2}- \sum
_{p=0}^{[\frac{\ell-1}{2}]}n_V(M({-\ell+2p})) -2\sum_{p=0
}^{[\frac{\ell}{2}]} n_V(T(\ell-2p)) - \sum_{p>0} n_V(T(\ell+2p)).
\end{align*} Using these formulas, we write down the non--zero
multiplicities in all cases. Below we assume $p\in \bz^+$.

\begin{enumerit}
\item  $V=V(r)\otimes M(s)$ with $s<0$.  $$n_V(T(r+s-2p)) = 1,\ \text{if}\ 2p\le r+s+1,\ \ n_V(M(r+s-2p)) = 1,\ \text{if} \
r+s+2\le p\le r.$$

\item $V=M(r)\otimes M(s)$ with $r<0$ or $s<0$.
$$n_V(T(r+s-2p)) = 1,\ \text{if}\ 2p\le r+s+1,\ \ n_V(M(r+s-2p))=1,\ \text{ if} \
p\ge r+s+2.$$

\item[(iii)] $V=V(r)\otimes T(s)$.
 $$n_V(T(r+s-2p)) =1, \ \text{if}\ p\le \min\{r,s\}, \ \
n_V(T(r+s-2p)) =2, \ \text{if}\ 2s+2\le 2p\le r+s+1.$$

\item[(iv)] $V=M(r)\otimes M(s)^{\sigma}$. $$n_V(T(r+s-2p)) = 1,  \ \text{if}\ 2p\le r+s+1,\ \
n_V(M(r+s-2p)) = 1,\ \text{if}\ p\ge r+s+2.$$

\item[(v)] $V=T(r)\otimes M(s)$.

$$n_V(T(r+s-2p)) = 1, \ \text{if}\ 2p\le \min\{2r,r+s+1\}, \ \ n_V(T(r+s-2p)) = 2, \ \text{if}\ 2r+2\le 2p\le r+s+1,$$
and
\begin{enumerit}
\item[(a)] $s\ge 0$:
 $$n_V(M(r+s-2p))= 1, \ \text{if}\ \max\{r,s\}+1\le p\le r+s+1, \ \ n_V(M(r+s-2p))= 2, \ \text{if}\ p\ge r+s+2,$$
\item[(b)] $s< 0$:
 $$n_V(M(r+s-2p))= 1, \ \text{if}\ r+s+2\le p\le r, \ \ n_V(M(r+s-2p))= 2, \ \text{if}\ p\ge r+1.$$
\end{enumerit}

\item[(vi)] $V=T(r)\otimes T(s)\cong T(r)\otimes M(s)
\oplus T(r)\otimes M(-s-2)$.

\end{enumerit}

\subsection{} Given $V\in \cal O$, let $B(V)$ be the branched crystal
$$B(V) = \text{\scriptsize $\bigoplus_W$ } B(W)^{\oplus n_V(W)},$$
where the sum is over the indecomposable objects $W\in \cal O$.
A straightforward comparison between sections \ref{tpdr} and \ref{tpdc}  shows the following.

\begin{thm} For all $U,V,W\in\cal O$ we have:
\begin{enumerit}
\item[(i)] $B(U)\otimes B(V)\cong B(U\otimes V)=B(V\otimes U)\cong B(V)\otimes B(U)$,
\item[(ii)] $B(U)\otimes
(B(V)\otimes B(W))\cong B(U\otimes (V\otimes W))=B((U\otimes
V)\otimes W)\cong (B(U)\otimes B(V))\otimes B(W)).$
\end{enumerit}
In particular,  the tensor product of branched crystals is
commutative and associative.
\end{thm}

\end{document}